\documentclass[reqno]{amsart}
\usepackage{amsmath, amsfonts, amsthm, amssymb}

\theoremstyle{plain}
\newtheorem{lemma}{Lemma}[section]
\newtheorem{theorem}{Theorem}[section]

\theoremstyle{definition}
\newtheorem{definition}{Definition}[section]
\theoremstyle{remark}
\newtheorem{remark}{Remark}[section]

\numberwithin{equation}{section}

\newcommand{\p}{\partial}

\newcommand{\di}{\mathrm{div}}
\newcommand{\grad}{\mathrm{grad}}

\newcommand{\fei}{\varphi}

\makeatletter
\newcommand{\rmnum}[1]{\mathrm{\romannumeral #1}}
\newcommand{\Rmnum}[1]{\mathrm{\expandafter\@slowromancap\romannumeral#1@}}
\makeatother

\begin{document}
\title[Uniqueness of Transonic Shocks in Nozzles]
{Global Uniqueness of Transonic Shocks in  Divergent Nozzles for
Steady Potential Flows}

\author{Li Liu}
\author{Hairong Yuan}
\address{Li Liu: Department of Mathematics, Shanghai Jiaotong
 University, Shanghai 200240, China}
\email{llbaihe@gmail.com}
\address{Hairong Yuan:
Department of Mathematics, East China Normal University, Shanghai
200241, China} \email{hryuan@math.ecnu.edu.cn,\,
hairongyuan0110@gmail.com}

\date{\today}

\begin{abstract}
We show that for steady compressible potential flow in a class of
straight divergent nozzles with arbitrary cross-section, if the flow
is supersonic and spherically symmetric at the entry, and the given
pressure (velocity) is appropriately large (small) and also
spherically symmetric at the exit, then there exists uniquely one
transonic shock in the nozzle. In addition, the shock-front and the
supersonic flow ahead of it, as well as the subsonic flow behind of
it, are all spherically symmetric.  This is a global uniqueness
result of free boundary problems of elliptic--hyperbolic mixed type
equations. The proof depends on the maximum principles and judicious
choices of comparison functions.
\end{abstract}
\maketitle


\section{Introduction and Main Results}

This paper is devoted to establishing uniqueness in the large for a
class of transonic potential flows with shocks in three--dimensional
divergent nozzles. These transonic shocks are spherically symmetric
and widely used in aerodynamics and computational fluid dynamics to
simulate those transonic shocks appeared in the so called de Laval
nozzles (i.e., convergent--divergent nozzles, cf.
\cite{an,CF,LY,Oe,Yu1}).

To formulate the problem, it would be somewhat convenient to use
$(r,\theta,\phi)$, the spherical coordinates of a point in
$\mathbf{R}^3$, with $0<\theta<\pi, 0<\phi<2\pi$. Let $\Sigma$ be a
$C^2$ domain on the unit $2$-sphere $\mathbf{S}^2\subset
\mathbf{R^3}.$ For two fixed positive constants $r^0<r^1$, we call
$\Omega=\{(r,\theta,\phi)\,|\, r^0< r< r^1, (\theta,\phi)\in
\Sigma\}$ a straight divergent nozzle, and $\Sigma$ its
cross--section. The wall of the nozzle is the truncated cone
$\Gamma=\{(r,\theta,\phi)\,|\, r^0< r< r^1, (\theta,\phi)\in
\p\Sigma\}$. For $i=0,1$, $\Sigma^i=\{(r,\theta,\phi)| r=r^i,
(\theta,\phi)\in \bar{\Sigma}\}$ is called the entry and the exit of
the nozzle respectively. Then obviously
$\p\Omega=\Gamma\cup\Sigma^0\cup\Sigma^1.$ Recall that the standard
Euclidean metric of $\mathbf{R}^3$ written in the local spherical
coordinates is $G=dr\otimes dr+r^2(d\theta\otimes
d\theta+\sin^2\theta d\phi\otimes d\phi).$

We consider steady potential polytropic gas flows in $\Omega$. The
governing equations are the following conservation of mass and
Bernoulli law (cf. \cite{CYu,CF}):
\begin{eqnarray}
&&\di(\rho\, \grad\fei)=0,\label{101}\\
&&\frac{1}{2}\langle\grad\fei,
\grad\fei\rangle+\frac{\rho^{\gamma-1}-1}{\gamma-1}=b_0. \label{102}
\end{eqnarray}
Since we use spherical coordinates rather than Descartesian
coordinates, here $\di,$ $\grad, \langle\cdot, \cdot\rangle$ are the
divergence operator, gradient operator and inner product with
respect to the metric $G$ respectively. The unknown $\fei$ is the
velocity potential (that is, $\grad\fei$ is the velocity of the
flow), $b_0$ is the known Bernoulli constant determined by the
incoming flow and/or boundary conditions, $\rho$ is the density, and
$\gamma>1$ is the adiabatic exponent. (Similar to \cite{CYu}, our
results and proof also hold for the isothermal case $\gamma=1$.) The
pressure of the flow $p$ and the speed of sound $c$ are determined
by
\begin{equation}\label{gamma-law}
p(\rho)=\frac{\rho^\gamma}{\gamma},\qquad c^2(\rho)=\rho^{\gamma-1}.
\end{equation}
Expressing $\rho$ in terms of $v=\sqrt{\langle\grad\varphi,
\grad\fei\rangle}=\sqrt{(\p_r\fei)^2+\frac{(\p_\theta\fei)^2}{r^2}+
\frac{(\p_\phi\fei)^2}{(r^2\sin^2\theta)}}$:
$$
\rho=\rho(v^2)
=\Big(1+(\gamma-1)(b_0-\frac{1}{2}v^2)\Big)^{\frac{1}{\gamma-1}},
$$
equation \eqref{101} becomes a second order equation for $\fei$:
\begin{equation}\label{103-a}
\di (\rho(\langle\grad\fei,\grad\fei\rangle)\ \grad\fei)=0.
\end{equation}
It is well known that this  equation is of mixed elliptic-hyperbolic
type in general \cite{CF}; it is elliptic if and only if the flow is
subsonic, i.e., $v<c$ or equivalently, by the Bernoulli law, $v<
c_*:=\sqrt{\frac{2}{\gamma+1}\big(1+(\gamma-1)b_0\big)}.$

Remember that we call $\Sigma^0, \Sigma^1$ the entry and exit of the
nozzle respectively.  This means that, we always assume:

\medskip
\noindent {\rm ($H$)} \qquad\qquad\qquad $\p_r\fei\ge0$ \qquad on
\,\,\, $\Sigma^i,\  i=0, 1$.

\medskip
\noindent That is, the gas flows in $\Omega$ on $\Sigma^0$ and flows
out of $\Omega$ on $\Sigma^1$.
\smallskip

Tremendous experiments and numerical simulations have shown that,
for a given supersonic flow near the entry of the nozzle, then by
giving an appropriately large back pressure at the exit, a transonic
shock must appear in the nozzle (see \S 147 in \cite{CF} or \S 4.3.4
in \cite{Oe}): the flow is discontinuous across a shock-front; the
flow ahead of the shock-front is supersonic, and behind of it is
subsonic, with pressure increases (velocity decreases) across the
shock-front. The position of the shock-front depends continuously
upon the back pressure and vice versa. Therefore, to design  nozzles
work for special purposes in aerodynamics, one has to understand the
existence, uniqueness, and stability of these transonic shocks via
rigorous theoretical analysis.

A basic  strategy to attack this problem is as follows: (a)
formulate a physically meaningful boundary value problem for the
governing partial differential equations (PDE), and look for some
special solutions involving transonic shocks; (b) study the
stability of these special solutions under perturbations of the
upcoming supersonic flow, the shape of the nozzle, and the back
pressure etc.; (c) study uniqueness of the special solutions in a
large class of functions. In step (a), instead of studying the full
Euler system, one always constructs some simplified models (such as
the potential flow equation, the popular quasi-one-dimensional model
of nozzle flow, see \cite{CF,Yu1,Yu4} and references therein), and
then might obtain solutions with particular symmetry by solving the
reduced ordinary differential equations or algebraic equations. Step
(b) essentially involves nonlinear small-perturbation problems: for
example, one  solves various linear PDE and then applying some
elegant nonlinear iteration techniques to show existence of the
nonlinear problem, see, for instance,
\cite{ChF1,ChF2,ChF3,CCF,CY,LY,XY,Yu2,Yu3}. In general, it is felt
that step (c) is more harder, which involves nonlinear problems
which are not of small-perturbation type.

This paper devotes exactly to establishing the uniqueness in the
large of a class of spherical symmetric transonic shocks in the
divergent nozzle $\Omega.$ Such spherical transonic shocks were
constructed in \cite{CF} and then widely used to explain the
transonic shock phenomena in de Laval nozzles (cf.
\cite{CF,Oe,Yu1}). Their stability, for the full steady Euler
system, in two-dimensional case, was studied by \cite{LY} recently.
See also \cite{XY1,XY2}. These results provide strong theoretical
supports for the applications of these special transonic shocks in
practices of aerodynamics, and in computational fluid dynamics to
test various numerical schemes designed to capture shocks (cf.
\cite{an,CF,Oe}).

We remark that in \cite{CYu}, Chen and Yuan have proved uniqueness
of a class of flat transonic shocks in straight ducts, and as a
byproduct, demonstrated that these flat shocks are instable,
therefore not physical. This result is also consistent with previous
instability results obtained in the papers \cite{CY,XY,Yu2}, which
were devoted to the stability issues by studying small-perturbation
problems.
\smallskip

Now let us formulate the boundary value problem of the potential
flow equation in $\Omega$ and summarize some important properties of
the special spherically symmetric transonic shocks. Then we will
state our main result. The proof is left in the next section.

We suppose that the flow is spherically symmetric and supersonic
(i.e., $v>c_*$) on $\Sigma^{0}$; spherically symmetric and subsonic
(i.e., $v<c_*$) on $\Sigma^1$. More specifically, for a constant
$u^0\in( c_*, \sqrt{2(b_0+\frac{1}{\gamma-1})})$ and a constant
$v_1\in(0, c_*)$, we consider the following problem:
\begin{eqnarray}
& \mbox{\eqref{103-a}}
                   &\text{in}\ \  \Omega,\label{103}\\
& \fei=u^0r^0, \quad  \p_{r}\fei=u^0   &\text{on}\ \  \Sigma^{0},
\label{104}\\
&\langle\grad\fei,\grad\fei\rangle=v_1^2  &\text{on}\ \  \Sigma^{1},
 \label{106}\\
&\langle\grad\fei, N\rangle=0 & \text{on}\ \  \Gamma, \label{0107}
\end{eqnarray}
where $N$ is the outward unit normal of $\Gamma$. (If the outward
unit normal of $\p\Sigma$ is $(n_1,n_2)$, then $N=(0, n_1/r^2,
n_2/r^2).$)

We remark that this problem is physically reasonable. In fact, since
our purpose is to study transonic shock phenomena in $\Omega$, the
flow should be supersonic at the entry, and the equation
\eqref{103-a} is hyperbolic there in the $r$-direction. So we need
two initial-value conditions like \eqref{104}. We choose $\fei$ as
in the first condition just to write the solution neatly. The flow
is supposed to be subsonic at the exit, where \eqref{103-a} is of
elliptic type, hence one and only one boundary condition is
necessary. Bernoulli type condition \eqref{106} means that the speed
of the flow is given at the exit, and by Bernoulli law, this is
equivalent to impose a uniform back pressure at the exit, which is
more physical than other conditions, such as $\fei$ itself (a
Dirichlet condition) (cf. \cite{CF}). In addition, \eqref{0107} is
the well known impenetrability or slip condition of inviscid flow
along solid boundary.

We will study uniqueness of solutions of \eqref{103}--\eqref{0107}
with the following structure (see \cite{CYu}).
\begin{definition}\label{d1}
For a $C^1$ function $r=f(\theta,\phi)$ defined on $\bar{\Sigma}$,
let
\begin{eqnarray*}
&& S=\{(f(\theta,\fei), \theta, \fei)\in \bar{\Omega} \,|\,
(\theta,\fei)\in \bar{\Sigma}\},\\
&&\Omega^-=\{(r,\theta,\fei)\in \Omega\, |\, r<f(\theta,\fei)\},\\
&&\Omega^+=\{(r,\theta,\fei)\in \Omega\,|\, r>f(\theta,\fei)\}.
\end{eqnarray*}
Then $\fei\in C^{0,1}(\bar{\Omega})\cap C^2(\overline{\Omega^-})\cap
C^{2}(\overline{\Omega^+})$ is a {\it transonic shock solution} of
\eqref{103}--\eqref{0107} if it is supersonic in $\Omega^-$ and
subsonic in $\Omega^+$, satisfies equation \eqref{103-a} in
$\Omega^-\cup \Omega^+$ and the boundary conditions
\eqref{104}--\eqref{0107} point-wise, the Rankine-Hugoniot jump
condition on $S$:
\begin{equation}
\rho(\langle\grad\fei^+,\grad\fei^+\rangle)\langle \grad\fei^+,
\nu\rangle=\rho(\langle\grad\fei^-,\grad\fei^-\rangle)
\langle\grad\fei^-, \nu\rangle,\label{rh}
\end{equation}
and the physical entropy condition on $S$:
\begin{equation}
\langle\grad\fei^+,\grad\fei^+\rangle<\langle\grad\fei^-,\grad\fei^-\rangle,
\label{entropy}
\end{equation}
where $\nu$ is the normal vector of $S$, and $\fei^+ (\fei^-)$ is
the limit value along $S$ of $\fei$ restricted in
$\overline{\Omega^+}$ ($\overline{\Omega^-}$). The surface $S$ is
also called the {\it shock-front}.
\end{definition}

For the existence of special transonic shock solutions  and their
important properties, we have the following lemma. Without any
ambiguity, we will also denote by $\fei^+ (\fei^-)$  the restriction
of $\fei$ in $\Omega^+ (\Omega^-).$

\begin{lemma}\label{l1}
Suppose the solution $\fei$ depends only on $r$. Then  for a given
$u^0$, we have the following results:

${\mathrm(1)}$ There is one solution $\fei^-(r)\in
C^2(\bar{\Omega})$ which is supersonic and solves problem
\eqref{103} \eqref{104} \eqref{0107} in the whole nozzle $\Omega$.
In addition, $\p_r\fei^-(r)$ is strictly monotonically increasing on
$(r^0,r^1).$

$\mathrm(2)$ There exists a connected open interval $I\subset
\mathbf{R}^+$ such that for $v_1\in I$, there is uniquely one
transonic shock solution
\begin{eqnarray}
\fei_b(r)=\begin{cases} \fei^-(r), & r\in[r^0,r_{s'}),\\
\fei_b^+(r), & r\in[r_{s'},r^1]
\end{cases}
\end{eqnarray}
to problem \eqref{103}--\eqref{0107}, with
$S_b:=\{(r,\theta,\phi)\in\Omega\, |\, r=r_{s'}\in(r^0,r^1)\}$ being
the shock-front, $\fei^-(r)$ the supersonic flow obtained in
$\mathrm(1)$, and $\fei_b^+(r)$ a subsonic flow. In addition, it has
the following properties:

$(\rmnum{1})$  $\p_r\fei_b^+(r)$ is strictly  monotonically
decreasing on $(r_{s'}, r^1)$.

$(\rmnum{2})$ $\fei^-(r_{s'})=\fei_b^+(r_{s'})$ and
$\fei^-(r)>\fei_b^+(r)$ for $r_{s'}<r \le r^1;$

$(\rmnum{3})$ $\p_r\fei^-(r)>\p_r\fei_b^+(r)$ for $r_{s'}\le r \le
r^1;$

$(\rmnum{4})$ For any fixed $R\in(r^0,r^1]$, $\p_r\fei_b^+(R)$ may
be regarded as a continuous function of $r_{s'}$, and is strictly
monotonically increasing for $r_{s'}\in(r^0,R).$
\end{lemma}

\begin{proof}
Let $v=\p_r\fei$ and solve $\rho_0$ from
$(u^0)^2/2+({\rho_0^{\gamma-1}-1})/({\gamma-1})=b_0.$ Consider the
following algebraic equations of $v(r)$ and $\rho(r)$ for
$r\in(r^0,r^1)$:
\begin{eqnarray}
&&r^2\rho v=a_0:=(r^0)^2\rho_0 u^0,\\
&&\frac 12v^2+\frac{\rho^{\gamma-1}-1}{\gamma-1}=b_0.
\end{eqnarray}
The claims in the Lemma can then be shown by elementary calculus
(see, for example, \cite{CF,Yu1}). In \cite{Yu1} the computation is
carried out for the two--dimensional steady full Euler system, but
the results there (i.e., Propositions 1, 3 and Theorem 6) still hold
for three-dimensional potential flows (except the asymptotic
expressions in Proposition 1 as $r\rightarrow\infty$, which we do
not need in this paper). Recall that $v=\p_r\fei$ here is the
velocity of the flow in the $r$-direction.
\end{proof}

We also have the global uniqueness of supersonic flow $\fei^-$ in
the whole nozzle $\Omega$, which may be proved by standard energy
estimates of nonlinear wave  equations.
\begin{lemma}\label{l2}
For a given $u^0$,  the solution $\fei^-(r)\in C^2(\bar{\Omega})$
which solves problem \eqref{103} \eqref{104} \eqref{0107} is unique
in the class of $C^2$ functions describing supersonic flows.
\end{lemma}

Now we state our main result.

\begin{theorem}\label{thm1}
Under the hypotheses $\mathrm(H)$, for a given $u^0$ and then any
$v_1\in I$, there exists one and only one transonic shock solution
to problem \eqref{103}--\eqref{0107} in the sense of Definition
\ref{d1}, which is exactly the $\fei_b(r)$ constructed in Lemma
\ref{l1} corresponding to $u^0$ and $v_1$.
\end{theorem}

\begin{remark}
This uniqueness result also holds if we consider the whole spherical
shell $\Omega'=\{(r,P)\, |\,
 r^0<r<r^1, P\in \mathbf{S}^2\}$ instead of the nozzle $\Omega$.
The proof is similar and simpler.
\end{remark}

Since the supersonic flow $\fei^-$ is known, independent of
downstream conditions and unique in the large according to Lemma
\ref{l2}, to prove Theorem \ref{thm1}, we just need to show that any
possible shock-front must be $S_b$ and the flow behind of it,
$\fei^+$, must be $\fei^+_b$. Therefore this is essentially a
uniqueness result of a free boundary problem for a nonlinear second
order PDE, with the transonic shock-front $S$ being the free
boundary. On the free boundary we have two boundary conditions,
namely the Dirichlet condition
\begin{eqnarray}\label{112}
&\fei^+=\fei^-, & \text{on}\ S
\end{eqnarray}
due to the continuity of $\fei$ across $S$ (see Definition
\ref{d1}), and the Neumann condition \eqref{rh}.

In the rest of this paper, Section 2, we will present the proof of
Theorem \ref{thm1}. It might be a little unexpected to find that the
proof is not very hard. It depends on maximum/comparison principles
of elliptic equations and some judicious choices of the above
constructed special solutions as comparison functions. The ideas are
generalizations of those in \cite{CYu}, and we need the deeper
properties of the family of special solutions as in Lemma \ref{l1}
to make them work. This fact, however, indicates that our methods
might be applicable to a large class of free boundary problems which
possess families of special solutions with fine structures.

\section{Proof of main result}

Let
\begin{eqnarray} \fei=\fei(r,\theta,\phi)=\begin{cases}
\fei^-(r), & r<f(\theta,\phi), \\
\fei^+(r,\theta,\phi), & r\ge f(\theta,\phi)
\end{cases}
\end{eqnarray}
be a transonic shock solution to problem \eqref{103}--\eqref{0107},
with $S=\{(r,\theta,\phi)\in\bar{\Omega}\, |\, r=f(\theta,\phi)\in
C^1(\bar{\Sigma}), (\theta,\phi)\in \bar{\Sigma}\}$ being the
shock-front. Then $\nu$, the normal of $S$, pointed from $\Omega^-$
to $\Omega^+$, in spherical coordinates, is
\begin{eqnarray}
\displaystyle \nu=\Big(1, -\frac{1}{r^2}\p_\theta f,
-\frac{1}{r^2\sin^2\theta}\p_\phi f\Big).
\end{eqnarray}
As in \cite{CYu}, by the Neumann condition \eqref{0107} and entropy
condition \eqref{entropy}, we can easily show that
\begin{eqnarray}
&\displaystyle\langle\nu, N\rangle=-\frac{1}{r^2}(n_1\p_\theta
f+n_2\p_\phi f)=0, & \text{at}\ \  {S}\cap\Gamma.
\end{eqnarray}
That is, the shock-front is always perpendicular to the wall of the
nozzle. In addition, direct computation yields
\begin{eqnarray}\label{204}
\rho\cdot \langle\grad\fei, \nu\rangle=\rho\cdot
\Big(\p_r\fei-\frac{1}{r^2}\p_\theta\fei\p_\theta
f-\frac{1}{r^2\sin^2\theta}\p_\phi\fei\p_\phi f\Big).
\end{eqnarray}

We may write the equation \eqref{103-a} in non-divergence form as
\begin{eqnarray}\label{205}
\sum_{i,j=0}^2 A^{ij}(D\fei)\p_{ij}\fei+B(D\fei)=0.
\end{eqnarray}
Here, for simplicity, we set $y^0=r, y^1=\theta, y^2=\phi$, and
$\p_{ij}=\p_{y^iy^j}.$ We do not need  to write out the specific
expressions of \eqref{205}, but it is essential to note that
$A^{ij}$ and $B$ depend only on the first order derivatives of
$\fei.$ This can be seen from the well known form of \eqref{205} in
Descartesian coordinates $(x^1,x^2,x^3)$:
\begin{eqnarray}
c^2\Delta\fei\label{pote}
-\sum_{i,j=1}^3\p_{x^i}\fei\p_{x^j}\fei\p_{x^ix^j}\fei=0.
\end{eqnarray}
Indeed, we can also use this expression below to show the validity
of maximum principles.

Now suppose for the given $u^0$ and $v_1\in I$, the special
transonic shock solution is $\fei_b$ and its shock-front is
$\{r=r_s\}\cap\bar{\Omega}.$ (See Lemma \ref{l1}.) The proof of
Theorem \ref{thm1} is then divided into two cases. \medskip

{\textsc{Case 1.}}\ We first show that if there holds
\begin{eqnarray}\label{206}
\min_{\bar{\Sigma}} f\ge r_s,
\end{eqnarray}
then $f\equiv r_s$ and
$w:=w(r,\theta,\fei)=\fei_b^+(r)-\fei^+(r,\theta,\phi)\equiv0.$

Let us consider the domain $\Omega_f^+=\{(r,\theta,\phi)\in\Omega\,
|\, r>f(\theta,\phi)\}$. By \eqref{206}, both $\fei^+$ and
$\fei_b^+$ are well defined in $\Omega_f^+.$ Therefore we may
formulate a boundary value problem of $w$ in $\Omega_f^+$ as
follows.

First, $w$ solves the following linear PDE:
\begin{eqnarray}\label{207}
&&\sum_{i,j=0}^2
a^{ij}(D\fei_b^+)\p_{ij}w+\sum_{i=0}^2b^i\p_iw\nonumber\\
&:=&\sum_{i,j=0}^2
A^{ij}(D\fei_b^+)\p_{ij}w\nonumber\\
&&+\Big(\sum_{i,j=0}^2\p_{ij}\fei^+(A^{ij}
(D\fei_b^+)-A^{ij}(D\fei^+))+B(D\fei_b^+)-B(D\fei^+)\Big)\nonumber\\
&=&0.
\end{eqnarray}
Note that there is no zeroth-order term $cw$ here. Since $\fei_b^+$
is subsonic, this equation is uniformly elliptic, and by our
requirements in Definition \ref{d1}, all  the coefficients $a^{ij},
b^i$ are bounded. So by the strong maximum principle \cite{GT}, if
$w$ is not a constant, its minimum can only be achieved on the
boundary
$\p\Omega_f^+=\Sigma^1\cup(\overline{\Omega_f^+}\cap\Gamma)\cup S.$

Second, by \eqref{106},  we have the following boundary condition of
$w$ on $\Sigma^1$:
\begin{eqnarray}\label{208}
\langle\grad\fei_b^++\grad\fei^+, \grad w\rangle=v_1^2-v_1^2=0.
\end{eqnarray}
By assumption $(H)$ and Lemma \ref{l1},
\begin{eqnarray}\label{ob}
\langle\grad\fei_b^++\grad\fei^+, (1,0,0)\rangle>0.
\end{eqnarray}
So \eqref{208} is a linear oblique derivative condition to $w$.
Similarly, on $\overline{\Omega_f^+}\cap\Gamma$, there is a Neumann
condition
\begin{eqnarray}\label{210}
\langle\grad w, N\rangle=0.
\end{eqnarray}
Therefore by Hopf boundary point lemma \cite{GT}, if $w$ is not
constant, the minimum of $w$ also can not be achieved on
$\Sigma^1\cup(\overline{\Omega_f^+}\cap\Gamma)$. (For points on
$\Sigma^1\cap(\overline{\Omega_f^+}\cap\Gamma)$, as in \cite{CYu},
we may use a locally even reflection arguments to prove that a
minimum can not be attained there, since the two surfaces meet there
at a right angle.)

Third, suppose the minimum is achieved at a point $(R,P)\in S$
(i.e., $R=f(P), P\in\bar{\Sigma}$). We note that $w$ satisfies a
Dirichlet condition here (cf. \eqref{112}):
\begin{eqnarray}\label{211}
w=g(\theta,\phi):=\fei_b^+(f(\theta,\phi))-\fei^-(f(\theta,\phi))\le
0.
\end{eqnarray}
The inequality holds due to property $(\rmnum{2})$ in Lemma
\ref{l1}, with $r_{s'}=r_s$.

We first suppose that $g$ admits a minimum at $P$, which is an
interior point of $\Sigma$. Then $\p_\theta g=\p_\phi g=0$ at $P$
and, since $\p_r\fei^-(r)>\p_r\fei_b^+(r)$ by property $(\rmnum{3})$
in Lemma \ref{l1} (taking $r_{s'}=r_s$), we also have
\begin{eqnarray}\label{grad0}
\p_\theta f(P)=\p_\phi f(P)=0,
\end{eqnarray}
which indicates that $\p_\theta w(R,P)=\p_\phi w(R,P)=0$, and
\begin{eqnarray}\label{qqq}
\p_\theta\fei^+(R,P)=\p_\phi\fei^+(R,P)=0.
\end{eqnarray} Indeed, by  \eqref{112} we have
$\fei^+(f(\theta,\phi),\theta,\phi)=\fei^-(f(\theta,\phi)),$ hence,
for example, we get  $(\p_r\fei^+)(\p_\theta f)+\p_\theta
\fei^+=\p_r\fei^-\p_\theta f$, and therefore $\p_\theta
\fei^+(R,P)=0$ by \eqref{grad0}.

Now  consider the Neumann condition on $S.$ By \eqref{rh},
\eqref{204}, \eqref{grad0} and \eqref{qqq}, there should hold
\begin{eqnarray}
\rho(|\p_r\fei^+|^2)\p_r\fei^+=\rho(|\p_r\fei^-|^2)\p_r\fei^-
\end{eqnarray}
at $(R,P)$. We may solve from this algebraic equation uniquely one
$\p_r\fei^+<\p_r\fei^-,$ which is not less  than $\p_r\fei_b^+$.
(Note that $R\ge r_s$. Here we used property $(\rmnum{4})$ in Lemma
\ref{l1}, with $r_{s'}=r_s$ and  $r_{s'}=R$.) Hence we have $\p_r
w=\p_r\fei_b^+-\p_r\fei^+\le0$ at $(R,P)\in S.$ However, we see that
$\nu=(1,0,0)$ pointed into $\Omega_f^+$ at $(R,P)$ by \eqref{grad0},
this is a contradiction to the Hopf boundary point lemma, which
asserts that there should hold $\p_rw>0$ at $(R,P)$, where $w$
attains its minimum.

For $P\in\p\Sigma,$ since $S$ is perpendicular to $\Gamma$, then
\eqref{grad0} still holds, and the above analysis also works.

Therefore $w$ must be a constant in $\overline{\Omega^+_f}.$ Now
look at \eqref{211}, we get $\p_\theta f=\p_\phi f\equiv 0$ for
$(\theta,\phi)\in \bar{\Sigma}.$ So $f$ is a constant. If $f=r_s,$
then $w=0$ and the uniqueness is proved.

Now we show that $f>r_s$ is impossible. Note that $w$ is a constant
and $\fei_b^+$ depends only on $r$, so $\fei^+$ depends only on $r$
for $f\le r\le r^1$, and therefore
\begin{eqnarray}
\fei(r,\theta,\phi)=\begin{cases} \fei^-(r), & r^0\le r<f,\\
\fei^+(r), & f\le r\le r^1
\end{cases}
\end{eqnarray}
is a special solution to problem \eqref{103}--\eqref{0107}. By Lemma
\ref{l1} $(\rmnum{4})$, since $f>r_s$, we must have
$v_1=\p_r\fei^+(r^1)>\p_r\fei_b^+(r^1)=v_1,$ a contradiction as
desired.

\medskip

{\textsc{Case 2.}}\ We now turn to the case that
\begin{eqnarray}\label{cases2}
R_s:=\min_{\bar{\Sigma}} f< r_s.
\end{eqnarray}
We will prove by contradiction that this is impossible.

Let
\begin{eqnarray}
\fei_B(r)=\begin{cases} \fei^-(r), & r^0\le r<R_s, \\
\fei_B^+(r), & R_s\le r\le r^1
\end{cases}
\end{eqnarray}
be a special transonic shock solution constructed in Lemma \ref{l1}
for which $\{r=R_s\}\cap \bar{\Omega}$ is the shock-front. Set
$V_1=\p_r\fei_B^+(r^1)$. Then by property $(\rmnum{4})$ in Lemma
\ref{l1} (with $R=r^1,$ $r_{s'}=R_s$ and then $r_{s'}=r_s$),
\begin{eqnarray}
V_1<v_1.
\end{eqnarray}

Consider the maximum of the function $w:=\fei_B^+-\fei^+$ defined in
$\Omega_f^+$. Now $w$ satisfies a linear uniformly elliptic equation
similar to \eqref{207}, only with $\fei_b^+$ replaced by $\fei_B^+.$
The oblique derivative condition on $\Sigma^1$ is now
\begin{eqnarray}\label{216}
\langle\grad\fei_B^++\grad\fei^+, \grad w\rangle=V_1^2-v_1^2<0,
\end{eqnarray}
and the Neumann condition on $\overline{\Omega_f^+}\cap\Gamma$ is
the same as \eqref{210}. So if $w$ is not a constant, its maximum
can only be attained on $S$, where the Dirichlet condition is
\begin{eqnarray}\label{217}
w=g(\theta,\phi):=\fei_B^+(f(\theta,\phi))-\fei^-(f(\theta,\phi))\le
0.
\end{eqnarray}
By the definition of $R_s$, the maximum, $0$, can be achieved at the
point $P\in\bar{\Sigma}$ where $f(P)=R_s$. Hence $\p_\theta
f(P)=\p_\phi f(P)=0$, as well as $\p_\theta \fei^+(R_s,P)=\p_\phi
\fei^+(R_s,P)=0$ and $\nu(R_s, P)=(1,0,0)$, as shown in Case 1.

By the Rankine-Hugoniot condition \eqref{rh}, $\p_r\fei^+$ should
satisfy
\begin{eqnarray}
\rho(|\p_r\fei^+|^2)\p_r\fei^+=\rho(|\p_r\fei^-|^2)\p_r\fei^-
\end{eqnarray}
at $(R_s,P)$. However, it follows from Lemma \ref{l1} that
$\p_r\fei_B^+(R_s)$ is the only solution to this algebraic equation
which satisfies the entropy condition. Therefore $\p_r w=0$ at
$(R_s,P)$, which is a contradiction to the Hopf boundary point lemma
from which  $\p_rw$ should be negative at $(R_s, P)$ if $w$ is not
constant.

Hence  $w$ is a constant and by \eqref{217}, it is zero.  Now we
obtain, from \eqref{216}, a contradiction. Therefore the Case 2 is
impossible.

This finishes the proof of Theorem \ref{thm1}.

\begin{remark}
In the proof of Case 1, we can also obtain contradiction if $w$ is
not constant by analyzing the maximum as in Case 2. But in Case 2 we
are restricted to considering only the maximum of $w$ to obtain
contradictions.
\end{remark}

\begin{remark}
The same proof also works for the uniqueness of cylindrical
transonic shocks for potential flows in two-dimensional straight
divergent nozzles.
\end{remark}

\begin{remark}
As noted in \cite{CYu}, in the proof we just used the fact from
Definition 1.1 that the flow on the right hand side of the
shock-front is subsonic (i.e., the entropy condition). We do not
need any assumption such as the flow should be subsonic in the whole
domain $\Omega_f^+.$
\end{remark}

\begin{remark}
The assumption $(H)$ is only used to guarantee the condition
\eqref{ob} at the exit to be oblique. Since for the family of
special solutions $\fei_b$ constructed in Lemma \ref{l1},
$\p_r\fei_b^+(r^1)$ is positive and bounded away from zero, say,
$\p_r\fei_b^+(r^1)>c_0>0$, therefore we can relax $(H)$ a little by,
for example, requiring only that $\p_r\fei\ge-c_0/2$ at the exit.

\end{remark}

\bigskip
{\bf Acknowledgments.} This research was supported in part by China
Postdoctoral Science Foundation $(20070410170)$, Shanghai Shuguang
Program (07SG29), Fok Ying Tung Foundation (111002), and the
National Science Foundation (USA) under Grant DMS-0720925. The
authors  thank sincerely  Professor Gui-Qiang Chen and Beixiang Fang
for their generous help and valuable comments.


\end{document}